\documentclass[letterpaper, 10 pt, conference]{ieeeconf} 
\IEEEoverridecommandlockouts 
\usepackage{cite}
\usepackage{amsmath,amssymb,amsfonts}
\usepackage{algorithmic}
\usepackage{graphicx}
\usepackage{algorithm,algorithmic}
\usepackage{hyperref}
\usepackage{textcomp}
\usepackage{xcolor}
\usepackage{tikz}
\usetikzlibrary{decorations.pathmorphing}

\newtheorem{theorem}{Theorem}
\newtheorem{assumption}{Assumption}
\newtheorem{definition}{Definition}

\newtheorem{remark}{Remark}

\begin{document}
\title{Dual MPC for Active Learning of Nonparametric Uncertainties}
\author{Tren Baltussen, Maurice Heemels, Alexander Katriniok
\thanks{The research has received funding from the European Union under the Horizon Europe Grant Agreement AIGGREGATE, no. 101202457.
All authors are with Control Systems Technology, Department of Mechanical Engineering, Eindhoven University of Technology, Eindhoven, The Netherlands. Contact: t.m.j.t.baltussen@tue.nl.\vspace{-0.6mm}}}

\maketitle

\thispagestyle{plain}
\pagestyle{plain}

\begin{abstract}
This manuscript presents a dual model predictive controller (MPC) that balances the two objectives of dual control, namely, system identification and control.
In particular, we propose a Gaussian process (GP)-based MPC that uses the posterior GP covariance for active learning.
The dual MPC can steer the system towards states with high covariance, or to the setpoint, thereby balancing system identification and control performance (exploration vs. exploitation).
We establish robust constraint satisfaction of the novel dual MPC through a contingency plan.
We demonstrate the dual MPC in a numerical study of a nonlinear system with nonparametric uncertainties.

Adaptive control, Optimal control
\end{abstract}

\section{Introduction}
One of the reasons for the high interest in model predictive control (MPC), in particular for safety-critical systems, is its ability to systematically handle constraints on the states and inputs.
Yet, the performance and safety of MPC heavily depend on the accuracy of the model of the system.
In various applications, the system dynamics may be uncertain or partially unknown, and need to be identified online.
This requires sufficiently informative data, and hence, the excitation of the system.
However, these system excitations may conflict with the control objective and safety requirements. Balancing these objectives, identification and control, is the focus of \textit{dual control} \cite{FELDBAUM1963541}. Therefore, studying dual control in the context of MPC is of high relevance for simultaneous identification and control of safety-critical systems.

\subsection{The Dual Effect}
When controlling a system that is subject to uncertainties, a controller can only make decisions based on the available information from the system.
However, the possibility of causally anticipating future measurements can be leveraged to improve state estimation or prediction, and hence, control performance.
When a controller incorporates a `future measurement program' and its associated statistics, this controller can causally anticipate these future measurements \cite{bar-shalom_dual_1974}.
The controller is said to have the dual effect, when, in addition to influencing the system state, it can also influence the state uncertainty, e.g., the state covariance.
As a consequence, the dual effect is comprised of two \mbox{properties \cite{bar-shalom_dual_1974}:}
\mbox{1) \textbf{Caution:} the} dependence of the state covariance on the control policy naturally introduces a degree of inherent caution -- or a lack thereof.
In the context of MPC, this dependence can be leveraged to tighten or relax the state constraints based on the anticipated state covariance.
\mbox{2) \textbf{Exploration:} this} anticipated state covariance can also be used to introduce an explorative or exciting element to the control policy.
Again, in the context of MPC, the state covariance can be used in the MPC cost function to steer towards regions of high uncertainty that may be of interest for system identification and control performance.
For a given control problem, the solution to the Bellman equation defines the globally optimal control policy, which naturally balances identification and control in an optimal way. In practice, however, solving the Bellman equation is generally intractable and requires approximate solutions \cite{mesbah_stochastic_2018}.

\subsection{Related Work}
\subsubsection{Dual Control}
Model predictive control is often implicitly used to approximate the Bellman equation in a tractable way by solving a finite-horizon optimal control problem \cite{rawlings2017model}.
In doing so, MPC relies on approximations of the true system dynamics and objective function. 
As a consequence, MPC may not fully capture the dual effect.
The dual effect can be incorporated in an approximative MPC policy either \textit{implicitly} or \textit{explicitly} \cite{mesbah_stochastic_2018}.
For a discussion and an overview of dual control in MPC, we refer the interested reader to \cite{mesbah_stochastic_2018}.

In implicit dual control, the MPC cost function is used to approximate the Bellman equation.
These types of methods are typically based on and related to approximate dynamic programming \cite{mesbah_stochastic_2018}.
For example, \cite{arcari_dual_2020a} and \cite{HUSHARP2024} use scenario trees in MPC to address the trade-off between exploration and exploitation. 
This exploration and exploitation trade-off is then \textit{implicitly} governed by the cost function \cite{mesbah_stochastic_2018}.
However, existing implicit dual MPC methods rely on additional safety filters \cite{HUSHARP2024}, are limited to linear systems \cite{soloperto_dual_2019}, or do not consider state constraints \cite{arcari_dual_2020a}, which is essential for safety-critical systems.
The implementation and computational complexity of implicit dual MPC limit its practical applicability \cite{mesbah_stochastic_2018}.

In explicit dual control, the explorative element of the control policy is introduced \textit{explicitly}, e.g., by enforcing a persistently exciting input signal, or adding a term to the cost function that minimizes (future) model uncertainty \cite{mesbah_stochastic_2018}.
We refer to \cite{mesbah_stochastic_2018} for a detailed review of these methods.
Although the exploration in explicit dual MPC typically requires some form of expert knowledge and is often heuristic \cite{mesbah_stochastic_2018}, these methods can explicitly account for state constraints. Therefore, explicit dual MPC is of high interest, especially for safety-critical systems.
While existing methods are highly relevant, they mostly focus on systems with parametric uncertainties \cite{mesbah_stochastic_2018} or fixed model structures \cite{heirung_stochastic_2017}.
This raises the question of how to effectively extend these methods to address nonparametric model uncertainties in dual MPC.

\subsubsection{Gaussian Process MPC}
Gaussian process (GP) regression is a Bayesian regression method that enables flexible, nonparametric modeling along with explicit uncertainty quantification, which can be leveraged with limited expert knowledge and few data points \cite{kocijan_modelling_2016}.
It has been successfully used to identify and control systems with unknown dynamics using MPC \cite{hewing_cautious_2020,baltussen2025,dubied2025robustadaptivempcformulation}.
While GP-MPC is often applied as an effective heuristic method \cite{hewing_cautious_2020,baltussen2025}, it can be complemented with robustly safe policies \cite{geurts_baltussen_2025} that ensure constraint satisfaction under uncertainty.
Alternatively, probabilisitic safety guarantees on GP-MPC are provided in \cite{dubied2025robustadaptivempcformulation}.
However, existing GP-MPC methods generally lack the dual effect as the GP covariance is not conditioned on the predicted control sequence.
Typically, the GP covariance is fixed prior to solving the MPC \cite{hewing_cautious_2020} or decoupled from the MPC in a zero-order optimization algorithm \cite{dubied2025robustadaptivempcformulation,lahr2024l4acados} to enable the real-time deployment of GP-MPC.
Although these works are of great interest, decoupling the GP covariance excludes the dual effect.
To bridge this gap, \cite{baltussen2025} presents a dual MPC that does condition the GP covariance on the predicted control sequence, yielding a controller with the dual effect, which is a neccessary property for active learning with GP-MPC.

\subsubsection{Active Learning}
Several works address active learning in MPC.
In \cite{soloperto_augmenting_2020}, an active learning framework augments existing MPCs while bounding deterioration of the primary cost. 
Similarly, \cite{soloperto_guaranteed_2023} adds a learning objective to the cost function with guarantees on closed-loop learning.
While these approaches promote exploration, they do not explicitly consider the dual effect, which is the focus of this work.

\subsection{Contributions}
As mentioned before, in \cite{baltussen2025} we introduced a GP-MPC that features the dual effect. However, we did not explicitly incentivize \textit{exploration}. Hence, this \textit{passively learning} dual MPC was limited to the \textit{cautious} property of the dual effect.
In this work, we extend the GP-MPC to \textit{actively} explore the state space and to learn unknown system dynamics while ensuring robust constraint satisfaction.
The main contributions of this manuscript are as follows.
1) We propose an explicit dual MPC that actively learns system dynamics using the active learning framework from \cite{soloperto_augmenting_2020}.
2)~We prove recursive feasibility of the actively learning MPC by incorporating a \textit{contingency horizon}\cite{alsterda_contingency_2019}.
3) We demonstrate the effectiveness of the proposed dual MPC through a numerical study on a nonlinear system with nonparametric, unknown dynamics.
In addition, we compare our dual MPC with a robust MPC (RMPC), a passively learning dual MPC \cite{geurts_baltussen_2025}, and the active learning framework from \cite{soloperto_augmenting_2020} with a single horizon.

\section{Problem Formulation}
\label{sec:Problem}
We consider a discrete-time control system $\mathcal{S}$ of the form
\begin{equation}
    \label{eq:system}
    x_{k+1} = f\left(x_{k}, u_{k} \right) + \underbrace{B g \left( x_{k}, u_{k} \right) +  v_k}_{=:w_k},
\end{equation}
where $x_k \in \mathbb{R}^{n_x}$ and $u_k \in \mathbb{R}^{n_u}$ denote the state and input vector, respectively, at time step $k \in \mathbb{N}$.
We assume to have access to exact measurements of the state $x_k$ and input $u_k$. \newline
The function $f: \mathbb{R}^{n_x} \times \mathbb{R}^{n_u} \rightarrow \mathbb{R}^{n_x}$ is a nominal model of the system dynamics.
The function $g: \mathbb{R}^{n_x} \times \mathbb{R}^{n_u} \rightarrow \mathbb{R}^{n_g}$ represents the unknown system dynamics that have to be identified.
We assume $f$ and $g$ are continuously differentiable and can be generally nonlinear.
The matrix $B \in \mathbb{R}^{n_x \times n_g}$ is a known matrix with full column rank. 
Lastly, $v_k \in \mathbb{V} \subset \mathbb{R}^{n_x}$ is a bounded disturbance at time $k \in \mathbb{N}$. The system is subject to state constraints $x_k\in \mathbb{X} \subset \mathbb{R}^{n_x}$ and input constraints $u_k \in \mathbb{U} \subseteq \mathbb{R}^{n_u}$, $k\in\mathbb{N}$. The uncertainty $w_k := B g(x_k,u_k)  + v_k$ is assumed to satisfy $w_k \in \mathbb{W} \subset \mathbb{R}^{n_x}$, for all $x \in \mathbb{X}$, $u \in \mathbb{U}$ and $k \in \mathbb{N}$, where $\mathbb{W}$ is a bounded set and $0\in\mathbb{W}$.

We consider the problem of designing a dual controller that balances system identification and control. Specifically, we aim to safely regulate the system $\mathcal{S}$ \eqref{eq:system} to a desired setpoint while actively learning the unknown dynamics $g(x,u)$.

\section{Preliminaries}
\subsection{Control Policies and the Dual Effect}
\label{sec:Dual_Control_Effect}
In the context of dual control, we can classify control policies based on the information that the controller utilizes \cite{bar-shalom_dual_1974}.
To this end, we introduce the concept of \textit{hyperstate}.
Let us denote the information available from system $\mathcal{S}$ at time $k \in \mathbb{N}$ by $\mathcal{I}_k := \{ x_k, \mathbf{z}_{k-1}, \dots, \mathbf{z}_{0} \}$, where $\mathbf{z}_i^\top = [x_i^\top, u_i^\top] \in \mathbb{R}^{n_z}$ are state–input pairs.
We define the hyperstate as the conditional probability of the state $x_k$ given the information $\mathcal{I}_k$ denote as $\xi_{k \mid k} := \mathrm{P} (x_k \mid \mathcal{I}_k)$.
Note that, since we assume exact state measurements, $\xi_{k \mid k}$ is a Dirac measure, but the general definition of the hyperstate allows for uncertainty in the state estimate, in the presence of measurement noise \cite{mesbah_stochastic_2018}.

In \textit{open-loop} control $u^{\mathrm{ol}}_k(\xi_{k \mid 0})$ the policy relies only on a-priori system information.
In \textit{feedback} control $u^{\mathrm{fb}}_k(\xi_{k \mid k})$ the policy relies on the information collected up to the current time $k$, but no knowledge of the future measurements is utilized by the controller.
In contrast, a dual control policy utilizes knowledge of future measurements to improve estimation or prediction, and control performance. For this, we define the causal anticipation of the hyperstate below.
\begin{definition}
    The \textit{causally anticipated hyperstate} at time $k + i \in \mathbb{N}$ based on the information $\mathcal{I}_k$ at time $k \in \mathbb{N}$ is
    \begin{equation}
        \vspace{-1mm}
        \xi_{i \mid i, \, k} := \mathbb{E} \left[ \xi_{k+i \mid k+i} \mid \xi_{k \mid k} ,  \mathcal{I}_k \right].
    \end{equation}
    Further, we denote $\xi_{k \mid_{k}^{k+N}} := ( \xi_{0 \mid 0, \, k}, \xi_{1 \mid 1, \, k}, \dots, \xi_{N \mid N, \, k} )$.
\end{definition} \vspace{0.5em}
Now, we can define the \textit{$N$-measurement feedback} control policy $u^{N\text{-}\mathrm{fb}}_k(\xi_{k \mid_{k}^{k+N}})$  that causally anticipates the future observations $\mathcal{I}_k$ to $\mathcal{I}_{k+N}$ from time step $k$ to $k+N$.
Hence, the $N$-measurement feedback (and closed-loop\footnote{In fixed end-time stochastic optimal control, \textit{closed-loop} policies anticipate all future measurements up to the end time $T$ \cite{bar-shalom_dual_1974}, whereas in MPC the horizon is typically shorter $N < T$.}) control policies affect not only the system state but also the predicted state uncertainty, i.e., the causally anticipated hyperstate \cite{bar-shalom_dual_1974,mesbah_stochastic_2018}.
As a consequence, these policies exhibit the dual effect, which we define next.
The dual effect was first formalized in \cite{bar-shalom_dual_1974}, and is defined in \cite{mesbah_stochastic_2018} as follows.

\begin{definition}
\label{def:DualControl}
    A control policy is said to have the \textit{dual effect} if it can affect, with non-zero probability, at least one $r^{\text{th}}$-order central moment of a state variable $(r \geq 2 )$.
\end{definition}
\vspace{0.5em}

According to Definition~\ref{def:DualControl}, only the $N$-measurement feedback and closed-loop control policies exhibit this dual effect \cite{bar-shalom_dual_1974,mesbah_stochastic_2018}.
This effect is a necessary property for active learning in GP-MPC, as it allows the MPC to \textit{actively} influence the posterior GP covariance and steer the system toward uncertain regions to learn the dynamics and enhance performance.

\subsection{Robust Model Predictive Control}
\label{sec:RMPC}
We use RMPC as a baseline for the dual MPC, as it provides a robust control policy that can be used to ensure safety and constraint satisfaction in the presence of uncertainty \cite{geurts_baltussen_2025}.
This general RMPC problem is formulated as follows:
\begin{subequations}
\label{eq:RMPC}
\begin{align}
    \min_{\bar{x}_{0 \mid k}, \bar{U}_k, }&~ J(x_k, \bar{x}_{0 \mid k}, \bar{U}_k) \\
    \hspace{-1.8mm}\text{s.t. }  \bar{x}_{i+1 \mid k} & = f(\bar{x}_{i \mid k},   \bar{u}_{i \mid k}), \ \,  i=0,1,\dots,N-1,  \\
    \bar{x}_{0 \mid k} & \in x_{k} \oplus S,\\
    \bar{u}_{i \mid k} &\in\bar{\mathbb{U}}_i, \qquad \qquad \hspace{2mm} i=0,1,\dots,N-1, \\
    \bar{x}_{i \mid k} & \in\bar{\mathbb{X}}_i, \qquad \qquad \hspace{2mm} i=0,1,\dots,N,
\end{align}
\end{subequations}
where, respectively,  $\bar{x}_{i \mid k}$ and $\bar{u}_{i \mid k}$ denote the $i^{\mathrm{th}}$ step-ahead prediction of the nominal state and input made at time $k \in \mathbb{N}$.
Furthermore, $\bar{U}_k := (\bar{u}_{0 \mid k}, \bar{u}_{1 \mid k}, \dots, \bar{u}_{N-1 \mid k},)$ denotes the nominal control sequence over the prediction horizon of length $N \in \mathbb{N}_{\geq 1}$.
The robust recursive feasibility of \eqref{eq:RMPC} typically relies on tightened constraint sets $\bar{\mathbb{X}}_i \subseteq \mathbb{X}$ and $\bar{\mathbb{U}}_i \subseteq \mathbb{U}$, and the use of a robust control invariant (RCI) set as the terminal state set $\mathbb{X}_N$ \cite{rawlings2017model}.
The resulting control input applied to $\mathcal{S}$ at time $k$ is $u_k(\xi_{k \mid k}) := \kappa(x_k, \bar{x}^*_{0 \mid k}) +\bar{u}_{0 \mid k}^{*}(\xi_{k \mid k})$, where an ancillary control law $\kappa$ and an RCI tube $S$ can be used as, e.g., in tube-based MPC \cite{rawlings2017model}.
We assume robust recursive feasibility to be satisfied for our RMPC in \eqref{eq:RMPC}:
\begin{assumption}
\label{assum:terminal_set}
The RMPC \eqref{eq:RMPC} is robustly recursively feasible such that, if there is a pair $(\bar{x}_{0 \mid k}, \bar{U}_k)$ that is feasible for \eqref{eq:RMPC}, for state $x_k$ at time $k \in \mathbb{N}$, then there is, for all $w_k\in\mathbb{W}$, a pair $(\bar{x}_{0 \mid k+1}, \bar{U}_{k+1})$ that is feasible for \eqref{eq:RMPC}, for the state $x_{k+1} = f(x_k, \kappa(x_k, \bar{x}_{0 \mid k}) +\bar{u}_{0 \mid k}) + w_k$ at time $k+1$.
\end{assumption}

\subsection{Gaussian Process Regression for MPC}
We use GP regression to predict the unknown residual dynamics $g(x,u)$ in \eqref{eq:system} through Bayesian inference. By imposing a Gaussian prior distribution on $g(x,u)$ we can obtain an tractable posterior distribution $\xi_{i, \mid i, k}$ that causally anticipates future measurements and is conditioned on the predicted control sequence, leading to a controller with the dual effect.
We briefly recall key concepts of GP regression; see \cite{rasmussen_gaussian_2006, kocijan_modelling_2016} for details on GPs and \cite{hewing_cautious_2020} for GP-MPC.

\subsubsection{Gaussian Process Regression}
Firstly, we take measurements of the residual dynamics $y_k \in \mathbb{R}^{n_g}$, where
\begin{equation}
    y_{k+1} = B^{\dag} \left( x_{k+1} - f(x_{k}, u_{k}) \right) = g\left(x_k, u_k \right) + \Tilde{v}_k \notag
\end{equation}
with $B^{\dag} = \left(B^{\top} B\right)^{-1}B^{\top}$ and $B^{\dag} v_k = \Tilde{v}_k \sim \mathcal{N}(0, \Sigma_v)$.
We use the training data $\mathbf{Z}_{k} = [\mathbf{z}_{0}, \dots, \mathbf{z}_{k-1}] \in \mathbb{R}^{n_z \times k}$ and
${\mathbf{y}_{k} = [y_{1}, \dots, y_{k} ] \in \mathbb{R}^{n_g \times k}}$ at time $k \in \mathbb{N}_{\geq 1}$ in order to approximate the residual dynamics by a GP $\tilde{\mathbf{d}} \hspace{-0.5mm} : \mathbb{R}^{n_z} \rightarrow \mathbb{R}^{n_g}$ such that
$g\left(x_k, u_k\right) \approx \tilde{\mathbf{d}}\left( \mathbf{z}_k \right)$.
We confine ourselves to scalar GPs $(n_g = 1)$ for the sake of readability.
Secondly, we impose a prior distribution on the GP through a universal kernel function $\mathfrak{K}\left(\mathbf{z}, \mathbf{z}'\right)$ \cite{kocijan_modelling_2016}. As we use the GP for model augmentation, we employ a zero mean prior.
Conditioning the prior GP on the information $\mathcal{I}_k$ yields the posterior GP.
Finally, we predict the residual dynamics over the MPC horizon by querying the posterior GP $\tilde{d}(\mathbf{z}) \in \mathbb{R}$
at a predicted test point $\mathbf{z}_{i \mid k} \in \mathbb{R}^{n_z}$ for $i = 0, 1,\dots, N-1$ at time $k$ with
\begin{equation}
    \mathrm{P} \bigl( \tilde{d} \mid \mathbf{z}, \mathcal{I} \bigr) = \mathcal{N} \left(d\left(\mathbf{z}, \mathcal{I}  \right), \Sigma^{d}\left(\mathbf{z}, \mathcal{I} \right) \right),
    \label{eq:dynamics_GP_mpc}%
\end{equation}
where $d$ and $\Sigma^{d}$ denote the predicted mean and covariance function of the posterior distribution of the GP, respectively,
\begin{subequations}
    \vspace{-1em}
\begin{align}
    d(\mathbf{z}) &= \mathbf{k}_{\mathbf{z Z}} \left(K_{\mathbf{ZZ}} + I \sigma^2_v \right)^{-1} \mathbf{y}, \label{eq:Sparse_Mean}\\
    \Sigma^d(\mathbf{z}) &= \mathfrak{K}_{\mathbf{z z}} - \mathbf{k}_{\mathbf{z Z}} \left(K_{\mathbf{ZZ}} + I \sigma^2_v \right)^{-1}  \mathbf{k}_{\mathbf{Z}\mathbf{z}} \label{eq:Sparse_Cov}.
\end{align}
\label{eq:Sparse_GP}%
\end{subequations}
Here $\mathfrak{K}_{\mathbf{z z}'} = \mathfrak{K}\left(\mathbf{z}, \mathbf{z}'\right)$, $\mathbf{k}_{\mathbf{Z z}} \in \mathbb{R}^{{k}}$ is the concatenation of the kernel function evaluated at the test point $\mathbf{z}$ and the training set $\mathbf{Z}$, where $[\mathbf{k}_{\mathbf{Z z}}]_i = \mathfrak{K}\left(\mathbf{z}_i,\mathbf{z}\right)$ and ${\mathbf{k}_{\mathbf{Z z}}^\top = \mathbf{k}_{\mathbf{z Z}}}$, and $K_\mathbf{Z Z} \in \mathbb{R}^{{k} \times {k}}$ is a Gram matrix, which satisfies ${[K_{\mathbf{Z Z'}}]_{ij} = \mathfrak{K}\left(\mathbf{z}_i,\mathbf{z}'_{j}\right)}$.
The Gaussian noise on $\mathbf{y} \in \mathbb{R}^{1 \times k}$ is typically handled by regularizing the Gram matrix with the noise variance $\sigma_v^2$, while the predictions of $d$ are noise free.
Propagating the GP \eqref{eq:dynamics_GP_mpc} over multiple steps is intractable as the nonlinear dynamics yield a non-Gaussian state distribution. Therefore, we use typical approximations to preserve the Gaussian distribution of the state, as described in \cite{hewing_cautious_2020,baltussen2025}. In addition, we employ a sparse GP \cite{snelson_sparse_2005} with $M$ inducing points equally spaced over the MPC horizon.

\begin{remark}
    Note that in Section~\ref{sec:Problem}, we assume bounded uncertainties $w$ and $v$ to ensure robust safety, whereas GP regression assumes Gaussian (and thus unbounded) uncertainties.
    Nevertheless, GP-MPC provides a meaningful heuristic by approximating the uncertainty set using a multiple of the GP’s posterior covariance \cite{soloperto_guaranteed_2023}.
    Given that the true system noise is bounded, this approximation allows us to retain robust safety guarantees while exploiting the GP’s nonparametric modeling and uncertainty quantification \cite{geurts_baltussen_2025}.
\end{remark}

\section{Passive Dual Model Predictive Control}
\label{sec:Dual_MPC}
\subsection{Gaussian Process Model Predictive Control}
We consider a \textit{passively learning} GP-MPC \cite{baltussen2025}, that we afterward extend with active learning. This passive GP-MPC uses the posterior GP mean and covariance to anticipate the future hyperstate and adapt the state constraints accordingly. The passive dual GP-MPC is formulated as follows:
\begin{subequations}
\label{eq:Prim_MPC}
\begin{align}
\min _{\hat{x}_{0 \mid k}, \hat{U}_k} \, J & \bigl(x_k, \hat{x}_{0 \mid k}, \hat{U}_k \bigr) \\
\hspace{-2.8mm} \text { s.t. } \hat{x}_{i+1 \mid k} &= f \left( \hat{x}_{i \mid k}, \hat{u}_{i \mid k} \right) + B d \left( \hat{\mathbf{z}}_{i \mid k}, \mathcal{I}_k \right), \notag \\
& \hspace{30.5mm} i = 0,1, \dots, N-1, \label{eq:MPC_Mean} \\
\Sigma_{i+1 \mid k}^x = & \begin{bmatrix} \nabla^\top \hspace{-1mm} f \left( \hat{x}_{i \mid k}, \hat{u}_{i \mid k} \right) & \hspace{-2.5mm} B \end{bmatrix} \Sigma_{i \mid k} \begin{bmatrix} \nabla^\top \hspace{-1mm} f\left( \hat{x}_{i \mid k}, \hat{u}_{i \mid k} \right) & \hspace{-2.5mm} B \end{bmatrix}^{\top}\hspace{-1mm}, \notag \\
& \hspace{30.5mm} i = 0,1, \dots, N-1, \label{eq:MPC_Cov} \\
\hat{x}_{0 \mid k} & \in x_{k} \oplus S,\\
\Sigma^x_{0 \mid k} &= \boldsymbol{0}, \\
u_{i \mid k} &\in \mathbb{U}, \qquad \qquad \qquad \; i = 0,1, \dots , N-1, \\
\hat{x}_{i \mid k} &\in \Hat{{\mathbb{X}}}_{i \mid k} \bigl(\hat{x}_{i \mid k}, \Sigma^x_{i \mid k} \bigr) , \hspace{2mm} i = 1,2, \dots , N.
\end{align}
\end{subequations}%
Here, $\Sigma_{i \mid k}$ denotes the joint covariance of the state $\Sigma^{x}_{i \mid k}$ and the GP $\Sigma^{d}_{i \mid k}$. For details we refer the reader to \cite{baltussen2025}.
We use the state covariance to heuristically adapt the constraints to $\Hat{\mathbb{X}}_{i \mid k} (\hat{x}_{i \mid k}, \Sigma^x_{i \mid k} ) := \{ \hat{x}_{i \mid k} \in \mathbb{X} \ominus \mathbb{S}(2 \sigma_{x, i \mid k}) \}$, where $\sigma_{x, i \mid k} \in \mathbb{R}^{n_x}$ is the predicted standard deviation of the state and $\mathbb{S}$ is a hyperrectangle.
Note that this MPC only leverages the \textit{cautious} property of dual control and does not actively explore the state space.
The posterior mean \eqref{eq:MPC_Mean} and covariance \eqref{eq:MPC_Cov} of the state define the \textit{causally anticipated hyperstate}:
$
    \xi_{i \mid i , k} := \mathcal{N} ( \hat{x}_{i \mid k}, \Sigma^{x}_{i \mid k} ), i = 1,2,\ldots,N,
$
where $\xi_{i \mid i , k}$ depends on the causal anticipation of the information that will be available up to time $k+i$ computed at time $k \in \mathbb{N}$ based on the information $\mathcal{I}_k$, see Section \ref{sec:Dual_Control_Effect}. 
As in \cite{baltussen2025}, we use a custom GP implementation in CasADi \cite{andersson_casadi_2019} that preserves the sensitivities of the GP covariance. 
We use CasADi’s symbolic framework to compute the hyperstate gradients.
Consequently, the GP-MPC conditions the covariance, and, hence, the hyperstate, on the control sequence $\hat{U}_k$, and can be solved using standard optimization solvers like IPOPT \cite{wachter_implementation_2006}.
The resulting policy features the dual effect.

\vspace{0.25em}
\begin{theorem}
The MPC policy from \eqref{eq:Prim_MPC} $u^{N\text{-}\mathrm{fb}}_k(\xi_{k \mid_{k}^{k+N}}) = \kappa(x_k, \hat{x}^*_{0 \mid k}) + \hat{u}^*_{0 \mid k}(\xi_{k \mid_{k}^{k+N}})$ is an $N$-measurement feedback policy that has the dual effect defined in Definition~\ref{def:DualControl}.
\end{theorem}
\begin{proof}
By definition, we have $\xi_{0 \mid 0 , k} := \xi_{k \mid k } = x_k$.
The $i$-step control sequence $\hat{U}^*_{0:i-1 \mid k}$ affects the causally anticipated hyperstate $\xi_{i \mid i, k}$ through the posterior mean $\hat{x}_{i \mid k}$ \eqref{eq:MPC_Mean} and posterior covariance $\Sigma^{x}_{i \mid k}$ \eqref{eq:MPC_Cov} of the predicted state for $i = 1, 2, \dots, N$ at time $k \in \mathbb{N}$.
Therefore, the control input $\hat{u}^*_{0 \mid k}$ depends on the  causally anticipated hyperstate $\xi_{k \mid_{k}^{k+N}}$ through \eqref{eq:Prim_MPC}.
Moreover, the control sequence $\hat{U}^*_{0:i-1 \mid k}$ affects the posterior covariance $\Sigma^{x}_{i \mid k}$ with non-zero probability such that, for $i = 1, 2, \dots, N$, we have
\vspace{-0.75mm}%
\begin{equation}
    \Sigma^{x}_{i \mid k} := \mathbb{E}\bigl[\Sigma_{k+ i }^{x} \mid \mathcal{I}_{k}, \hat{U}^*_{0:i-1 \mid k} \bigr] \neq \mathbb{E}\bigl[ \Sigma_{ k+ i}^{x} \mid \mathcal{I}_{k} \bigr]. \notag
\end{equation}
\vspace{-0.75mm}%
Hence, $u^{N\text{-}\mathrm{fb}}_k(\xi_{k \mid_{k}^{k+N}})$ has the dual effect.
\end{proof}

\subsection{Safety Guarantees}
\label{sec:CMPC}%
We ensure robust constraint satisfaction via a contingency MPC as in \cite{geurts_baltussen_2025} and consider the following problem:
\begin{subequations}
\begin{align}
  \min_{\bar{x}_{0 \mid k}, \bar{U}_k,\hat{U}_k}&~ J_C(x_k, \bar{x}_{0 \mid k}, \bar{U}_k,\hat{U}_k)\label{subeq:cost}\\
  \hspace{-1.8mm}\text{s.t. }  \bar{x}_{i+1 \mid k} & = f(\bar{x}_{i \mid k},\bar{u}_{i \mid k}), \hspace{7.8mm} i=0,1,\dots,N-1, \label{subeq:fr} \\
  \hat{x}_{i+1 \mid k} & = f(\hat{x}_{i \mid k},\hat{u}_{i \mid k}) + B d(\hat{\mathbf{z}}_{i \mid k},\mathcal{I}_k) \notag \\  
  & \hspace{31.6mm} i=0,1,\dots,N-1,\label{subeq:fL}  \\
  \Sigma_{i+1 \mid k}^x = & \begin{bmatrix} \nabla^\top \hspace{-1mm} f ( \hat{x}_{i \mid k}, \hat{u}_{i \mid k} ) & \hspace{-2.5mm} B \end{bmatrix} \Sigma_{i \mid k} \begin{bmatrix} \nabla^\top \hspace{-1mm} f( \hat{x}_{i \mid k}, \hat{u}_{i \mid k} ) & \hspace{-2.5mm} B \end{bmatrix}^{\top}\hspace{-1mm}, \notag \\
  & \hspace{31.6mm} i = 0,1, \dots, N-1, \\
  \bar{x}_{0 \mid k} & = \hat{x}_{0 \mid k} \in x_{k} \oplus S
  , \label{subeq:same_initial_x}\\
  \Sigma^x_{0 \mid k} &= \boldsymbol{0}, \\
  \bar{u}_{0 \mid k} & = \hat{u}_{0 \mid k}, \label{subeq:same_initial_u}\\
    \bar{u}_{i \mid k} &\in\bar{\mathbb{U}}_i, \qquad \qquad \qquad \; i=0,1,\dots,N-1, \label{subeq:input_r} \\
    \hat{u}_{i \mid k} &\in\bar{\mathbb{U}}_i, \qquad \qquad \qquad \; i=0,1,\dots,N-1, \label{subeq:input_l} \\
   \bar{x}_{i \mid k} & \in\bar{\mathbb{X}}_i, \qquad \qquad \qquad \; i=0,1,\dots,N, \label{subeq:state_r} \\
   \hat{x}_{i \mid k} & \in\hat{\mathbb{X}}_{i \mid k} \bigl(\hat{x}_{i \mid k}, \Sigma^x_{i \mid k} \bigr), \hspace{3.2mm} i=0,1,\dots,N, \label{subeq:state_l}
\end{align}
\label{eq:MPC_problem}%
\end{subequations}
where the contingency MPC cost function is defined as:
\begin{equation}
\begin{aligned}
\label{eq:cost_contingency}
    J_C(x_k, \bar{x}_{0 \mid k}, \bar{U}_k,\hat{U}_k) = (1-\lambda) \, & J(x_k, \hat{x}_{0 \mid k}, \hat{U}_k) \, + \\  \lambda \, & J(x_k, \bar{x}_{0 \mid k}, \bar{U}_k),
\end{aligned}
\end{equation}
where $\lambda \in [0,1]$ is a parameter that weights the \textit{contingency horizon} \eqref{subeq:fr} and the \textit{performance horizon} \eqref{subeq:fL}.
This MPC uses a secondary robust input sequence $\bar{U}_k$, based on the RMPC from \eqref{eq:RMPC}, to ensure the safety of the GP-MPC.
By learning the dynamics $g$, the performance horizon can select a `better' control input $\hat{u}^*_{0 \mid k}$, while the contingency horizon restricts the input $\hat{u}^*_{0 \mid k} = \bar{u}^*_{0 \mid k}$ to the robustly admissible input set.
Consequently, the contingency MPC can reduce the conservatism of the RMPC while ensuring safety \cite{geurts_baltussen_2025}.

\begin{assumption}
\label{assum:disturbance}
Let $\mathbb{D}_{i \mid k} := \{ d(\mathbf{z}_{i \mid k},\mathcal{I}_k) \mid \mathbf{z}_{i \mid k} \in \mathbb{X} \times \mathbb{U} \}$ denote the set of all possible GP mean predictions at prediction step $i$ given the data $\mathcal{I}_k$ and $\hat{\mathbb{X}}_{i \mid k} := \{ \hat{x}_{i \mid k} \in \mathbb{X} \ominus 2 \sigma_{x, i \mid k} \}$ the tightened state constraints of the performance horizon.
The learned residual dynamics satisfy $\mathbb{D}_{i \mid k} \subseteq \mathbb{W}$, and it holds that $\hat{\mathbb{X}}_{i \mid k} \supseteq \bar{\mathbb{X}}_i$, for $i = 0,1,\dots,N$ and for all $k\in\mathbb{N}$.
\end{assumption}

\vspace{0.5em}
\begin{theorem}
\label{thm:RF_CMPC}
Let Assumptions \ref{assum:terminal_set} and \ref{assum:disturbance} hold. 
If OCP \eqref{eq:MPC_problem} has a feasible solution $(\bar{x}_{0 \mid k},\bar{U}_{k}, \hat{U}_{k})$ for $x_k$ at time $k\in\mathbb{N}$, then \eqref{eq:MPC_problem} is feasible for the next state $x_{k+1}=f(x_k, \kappa(x_k, \bar{x}_{0 \mid k}) +\bar{u}_{0 \mid k})+w_k$, for any $ w_k\in\mathbb{W}$, at time $k+1$.
\end{theorem}
\begin{proof}
    See \cite{geurts_baltussen_2025}.
\end{proof}

\vspace{0.5em}
\begin{remark}
    Note that if Assumption \ref{assum:disturbance} does not hold, soft-constraints can be employed in the \textit{performance horizon} to ensure that the contingency MPC remains feasible without sacrificing the safety guarantees \cite[Theorem 2]{geurts_baltussen_2025}.
\end{remark}

\section{Active Dual Model Predictive Control}
\subsection{Active Learning Framework}
We extend the passive GP-MPC \eqref{eq:MPC_problem} to leverage the posterior GP covariance for the active learning of nonparametric uncertainties.
To this end, we extend the active learning framework from \cite{soloperto_augmenting_2020} as it provides an intuitive framework to augment existing MPCs with active learning.
By combining this with the multi-horizon contingency MPC, we can leverage the dual effect of the performance horizon to actively learn the model while the contingency horizon ensures safety.
This framework uses a learning cost function $H$ while the primary cost function $J$ appears in additional constraints.

\vspace{-0.5em}
\begin{subequations}
\label{eq:AL_MPC}
\begin{align}
\min _{\hat{x}_{0 \mid k}, \bar{U}_k, \hat{U}_k, \Delta_k} \, & H(x_k, \hat{x}_{0 \mid k}, \hat{U}_k ) \\
\text { s.t. } & J (x_k, \hat{x}_{0 \mid k}, \hat{U}_k ) = J_k^B + \Delta_k,\\
& \Delta_k \leq \bar{\beta} \max \{ J^+_k, 0 \} + \Bar{\gamma} + Y_{k-1}, \label{eq:lim_det_mean}\\
& \Delta_k \leq \beta^{\text{max}} \max \{J^+_k, 0 \} + \gamma^{\text{max}}, \label{eq:lim_det_max}\\
& \eqref{subeq:fr} - \eqref{subeq:state_l}.
\end{align}
\end{subequations}
Here, the primary cost function $J$ of the performance horizon is composed of a baseline cost $J_k^B$, and a deterioration term $\Delta_k$ which quantifies the cost incurred for exploration. The deterioration used for exploration $(\Delta_k)$ is bounded over time by \eqref{eq:lim_det_mean} and at each time step by \eqref{eq:lim_det_max}. 
A storage variable $Y_{k} := Y_{k-1} + \bar{\beta} \, \max(J^+_k, 0) + \bar{\gamma} - \Delta_k$ accumulates cost deterioration, preventing excessive performance loss over time. The term $J^+_k := J(x_{k-1}, \tilde{x}_{0 \mid k - 1}, \tilde{U}_{k-1}) - J^B_k$, where $\tilde{x}_{0 \mid k - 1}, \tilde{U}_{k-1}$ denote the previous solution of \eqref{eq:AL_MPC}, allows exploration when it can potentially improve performance \cite{soloperto_augmenting_2020}. The learning behavior is tuned via the user-defined parameters $\bar{\beta}$, $\bar{\gamma}$, $\beta^{\text{max}}$, and $\gamma^{\text{max}}$, see \cite{soloperto_augmenting_2020} for details.
We define the baseline cost $J_k^B$ as the optimal cost $J(x_k, \hat{x}^*_{0 \mid k}, \hat{U}^*_{0 \mid k})$ from the performance horizon of \eqref{eq:MPC_problem}, which requires solving both \eqref{eq:MPC_problem} and \eqref{eq:AL_MPC} at each time step.
Alternatively, an upper bound on $J^B_k$ can be computed in different ways (see \cite{soloperto_augmenting_2020}).

In the GP-MPC framework, this active learning formulation steers the system to states that maximize model uncertainty with the aim to learn the unknown dynamics. To this end, the learning cost $H$ is defined as the negative sum of the posterior GP covariance over the prediction horizon:
\begin{equation}
\label{eq:Learn_Obj}
    H(x_k, \hat{x}_{0 \mid k}, \hat{U}_k ) = -\sum_{i=0}^{N} \Sigma^{d}( \hat{\mathbf{z}}_{i \mid k}, \mathcal{I}_k ).
\end{equation}
The resulting $N$-measurement feedback-dual control policy
\mbox{$u^{N\text{-}\mathrm{fb}}_k(\xi_{k \mid k + N} ) := \kappa(x_k, \hat{x}^*_{0 \mid k}) + \hat{u}^*_{0 \mid k}$} with $\hat{u}^*_{0 \mid k}$ obtained from \eqref{eq:AL_MPC} is used to balance system identification and control performance of the closed-loop system $\mathcal{S}$.

\subsection{Safety Guarantees}
Similar to \cite{soloperto_augmenting_2020}, the robust recursive feasiblity of the actively learning dual MPC is ensured by relying on a robust MPC formulation. In contrast to \cite{soloperto_augmenting_2020}, we only require constraints to be robustly enforced for the contingency horizon \cite{geurts_baltussen_2025}, which can reduce conservatism.
The robust recursive feasibility of the contingency MPC \eqref{eq:MPC_problem} is then ensured by using a robust MPC in the \textit{contingency horizon} \cite{geurts_baltussen_2025}.

\vspace{0.5em}
\begin{theorem}
\label{thm:RF_Active_CMPC}
Let Assumptions \ref{assum:terminal_set} and \ref{assum:disturbance} hold. 
If OCP \eqref{eq:AL_MPC} has a feasible solution $(\bar{x}_{0 \mid k},\bar{U}_{k}, \hat{U}_{k}, \Delta_{J_k})$ for $x_k$ at time $k\in\mathbb{N}$, then \eqref{eq:AL_MPC} is feasible for the next state $x_{k+1}=f(x_k, \kappa(x_k, \bar{x}_{0 \mid k}) +\bar{u}_{0 \mid k})+w_k$, for any $ w_k\in\mathbb{W}$, at time $k+1$.
\end{theorem}

\begin{proof}
    If the problem \eqref{eq:MPC_problem} has a feasible solution $(\bar{x}_{0 \mid k}, \bar{U}_k, \hat{U}_k)$ at time $k \in \mathbb{N}$, it is trivial that \eqref{eq:AL_MPC} has an feasible solution $(\bar{x}_{0 \mid k}, \bar{U}_k, \hat{U}_k, \Delta_k)$ with $\Delta_k = 0$.
    It remains to show, that the primary MPC \eqref{eq:MPC_problem} is feasible for next state $x_{k+1}=f(x_k, \kappa(x_k, \bar{x}_{0 \mid k}) +\bar{u}_{0 \mid k})+w_k$, for any $ w_k\in\mathbb{W}$, at the next time $k+1$. 
    Since every trajectory satisfying \eqref{eq:AL_MPC} must satisfy the constraints of \eqref{eq:MPC_problem} at time $k \in \mathbb{N}$, it follows from Theorem \ref{thm:RF_CMPC} that there exists a feasible solution $(\bar{x}_{0 \mid k+1}, \bar{U}_{k+1}, \hat{U}_{k+1})$ to \eqref{eq:MPC_problem} at time $k+1$.
    \end{proof}

\section{Numerical Example}
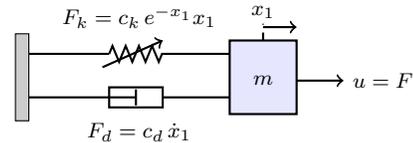
\begin{figure}[t]
    \centering
    \begin{tikzpicture}[scale=0.89, every node/.style={font=\footnotesize}]
        \draw[fill=gray!40] (-2.2,0) rectangle (-2,1.3);
        \draw[fill=blue!10,thick] (1,0.1) rectangle (2,1.2);
        \node at (1.5,0.6) {$m$};
        \draw[thick] (-2,1.0) -- (-0.8,1.0);
        \draw[thick,decorate,decoration={zigzag,segment length=5,amplitude=3}] (-0.8,1.0) -- (0,1.0);
        \draw[thick] (0,1.0) -- (1,1.0);
        \draw[->,thick,black] (-0.9,0.8) -- (0,1.2);
        \draw[thick] (-2,0.35) -- (-0.4,0.35);
        \draw[thick] (-0.8,0.5) rectangle (0,0.2);
        \draw[thick] (-0.4,0.23) -- (-0.4,0.47);
        \draw[thick] (0,0.35) -- (1,0.35);
        \node[below] at (-0.35,0.1) {$F_d = c_d \, \dot{x}_1$};
        \node[above] at (-0.35,1.25) {$F_k = c_k \, e^{-x_1} x_1$};
        \draw[dashed,thick] (1.5,1.2) -- (1.5,1.45);
        \node[above] at (1.5,1.4) {$x_1$};
        \draw[->,thick,black] (1.5, 1.4) -- (2, 1.4);
        \draw[->,thick,black] (2,0.6) -- (2.7,0.6);
        \node[right,black] at (2.7,0.6) {$ u = F$};
    \end{tikzpicture}
    \vspace{-0.5em}
    \caption{Schematic of the mass-spring-damper system.}
    \label{fig:system_schematic}
    \vspace{-1.5em}
\end{figure}

We perform a numerical study on a nonlinear mechanical system to compare the discussed MPC algorithms and demonstrate the proposed actively learning dual MPC. Specifically, we consider a mass–spring–damper system with a nonparametric unknown spring force, shown in Fig.~\ref{fig:system_schematic}, that is described by the continuous-time state-space equation:
\begin{equation}
    \label{eq:example_system}
    \begin{bmatrix}
        \dot{x}_1 \\ \dot{x}_2
    \end{bmatrix} =
    \underbrace{\begin{bmatrix}
        x_2 \\ - \frac{1}{m} \left( c_d \, x_2 + c_k \, x_1 + u \right)
    \end{bmatrix}}_{= \, f_c(x,u)} + \underbrace{\begin{bmatrix}
        0 \\ - \frac{1}{m}  c_k (e^{-x_1} - 1) x_1
    \end{bmatrix}}_{= \, B g_c(x,u)}, \notag
\end{equation}
where $x_1$ [m] is the displacement, $x_2$ [m/s] is the velocity, $m = 1$ [kg] is the mass, $c_d = 1$ [Ns/m] is the damping constant and $c_k = 0.33$ [N/m] is the spring constant.
We consider $v_k=0$.
The system is subject to state constraints $\mathbb{X} = [-0.1, 1.1] \times [-5, 5]$, and input constraints $\mathbb{U} = [-5, 5]$.
The goal is to control the system to a desired setpoint $x_{\text{ref}}^\top = \begin{bmatrix}1 & 0\end{bmatrix}$ and $x_{\text{ref}}^\top = \begin{bmatrix}1.095 & 0\end{bmatrix}$ after $t = 5$ s, while actively learning the unknown nonlinear spring force using the proposed dual MPC.
The continuous-time system is linearized at the origin and discretized using the forward Euler method with a sampling time of $0.1$ seconds, yielding the nominal discrete-time model $f$.
The input is applied with a zero-order hold and is kept constant over the sampling period.
The primary cost function is a standard tracking cost
\begin{equation}
    \vspace{-1mm}
    \label{eq:cost_primary}
    J(x_k, x_{0 \mid k}, U_k) = \sum_{i=0}^{N-1} \ell(x_{i \mid k}, u_{i \mid k}) + \ell_N(x_{N \mid k}),
    \vspace{-0.5mm}
\end{equation}
with stage cost $\ell(x, u) = \|x - x_{\text{ref}}\|_Q^2 + \|u \|_R^2$ and terminal cost 
$\ell_N(x) = \|x - x_{\text{ref}}\|_P^2$, and $Q \succeq 0$ and $R,P \succ 0$.

\subsection{MPC Implementation and Properties}
We consider the primary cost function \eqref{eq:cost_primary} with weighting matrices $Q = \text{diag}(20, 1)$ and $R = 1$ with $P$ from the LQR corresponding to $Q$ and $R$, and $\lambda = 10^{-3}$.
We consider a prediction horizon of $N = 20$.
For the GP regression, we use $x_1$ as the regression feature with a squared exponential kernel $\mathfrak{K}$ with a signal variance $\sigma_f^2 = 0.33$, noise variance $\sigma_v^2 = 0$ and length scale $l = 0.3$ \cite{kocijan_modelling_2016}. We use $M=4$ inducing points.
The unmodeled spring force only affects the velocity dynamics, hence, $B = [0, 1]^\top$.
We consider $\kappa \equiv 0$ and $w_k \in \mathbb{W}$ for all $k \in \mathbb{N}$ with $\mathbb{W} = \{0\} \times [-0.33, 0.33]$.
The parameters of the active learning framework are set to $\bar{\beta} = 1$, $\bar{\gamma} = 0$, $\beta^{\text{max}} = 0$, and $\gamma^{\text{max}} = \infty$, so that the degree of excitation decays over time as the tracking cost decreases.
The MPC parameters are constant for all MPC formulations.

The RMPC \eqref{eq:RMPC} robustly accounts for all possible realizations of the unmodeled dynamics and enforces robustly tightened constraints with an RCI terminal set, computed using the LQR, to ensure recursive feasibility.
The GP-MPCs leverage the online collected measurements to identify and predict the unmodeled dynamics and use adaptive constraint tightening based on the anticipated state covariance.
Furthermore, the GP-MPCs use the RMPC \eqref{eq:RMPC} in the contingency horizon.
The passive GP-MPC \cite{geurts_baltussen_2025} only leverages the cautious property of the dual effect.
The proposed dual GP-MPC \eqref{eq:AL_MPC} uses the learning objective \eqref{eq:Learn_Obj} to leverage both the cautious and explorative properties of the dual effect (Section \ref{sec:Dual_Control_Effect}).
Finally, we evaluate an implementation of the active learning framework by \cite{soloperto_augmenting_2020} that does not use a secondary horizon and relies on robustly tightened constraints from \eqref{eq:RMPC}, and, therefore, lacks the cautious property of the dual effect.

\subsection{Results}
\begin{figure}[t]
    \centering
    \includegraphics[width=1\linewidth]{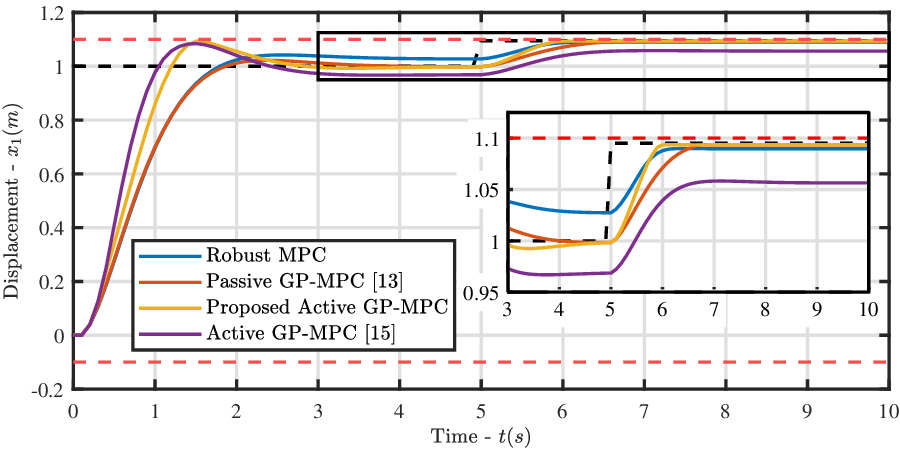}
    \vspace{-1.9em}
    \caption{Results of simulation study. The black and red dashed line indicate the MPC setpoint and state constraints, respectively. The setpoint changes after 5 seconds and the MPC does not have a preview of this change.}
    \label{fig:results}
    \vspace{-1.5em}
\end{figure}

The results of the simulation study are shown in \mbox{Fig. \ref{fig:results}}.
The RMPC considers the worst-case disturbance caused by the uncertain spring force.
Initially, the RMPC overshoots and has a large steady-state error of 2.74\%. After \mbox{$t=5$~s,} the improved tracking is largely incidental.
As the nominal model is unaware of the reduced spring stiffness, the controller tends to overcompensate, and the system happens to settle near the constraint, yielding a small steady-state error of 0.48\%, at $t = 10$~s.
The actively learning GP-MPC with a single horizon \cite{soloperto_augmenting_2020} actively explores the state space and later exploits the learned model.
Although the learned model is aware of the reduced stiffness, it is still constrained to the robust constraints that were designed for the worst-case disturbance. Consequently, the MPC cannot fully exploit the learned model, leading to conservative behavior and steady-state errors of 3.16\% and 3.53\%, at $t=5,10$~s, respectively.

The passively learning GP-MPC with the contingency horizon \cite{geurts_baltussen_2025} is able to successfully control the system to the desired setpoint while learning and adapting to the nonparametric unmodeled dynamics.
As the two horizons branch after the first input, the contingency horizon ensures that the RCI set is reachable at all times, while the performance horizon can exploit the learned model to reduce conservatism and improve performance, significantly reducing the steady-state errors to 0.17\% and 0.13\%, at $t=5,10$~s, respectively.
Lastly, the proposed actively learning GP-MPC can safely excite the system to learn the unmodeled dynamics, and exploit the learned model to obtain similar performance (0.19\% and 0.13\%).
The active GP-MPC \eqref{eq:AL_MPC} has an average solve time of $0.295$ seconds on a standard laptop PC using MATLAB, CasADi \cite{andersson_casadi_2019} and IPOPT \cite{wachter_implementation_2006}, demonstrating the practical computational efficiency of the proposed approach.

\section{Conclusion}
In this work, we proposed a Gaussian process (GP)-based dual MPC that leverages the dual effect to actively learn nonparametric uncertainties.
We extended the active learning framework of \cite{soloperto_augmenting_2020} by including the dual effect and adopted a multi-horizon contingency MPC formulation \cite{geurts_baltussen_2025} to ensure robust constraint satisfaction.
The proposed dual MPC features both caution and exploration, the two properties of dual control, enabling safe active learning and control near constraints, and it was demonstrated on a nonlinear system.

Although a simple problem was considered for the sake of analysis, 
we expect that the actively learning dual MPC will outperform the passive variant in more complex problems that require exploration to learn the unknown dynamics.
Furthermore, the results show that the contingency framework can reduce conservatism compared to a single-horizon MPC.
Future work will focus on the performance analysis of both the actively and passively learning dual MPC and explore the use of robust adaptive MPC to further reduce conservatism.

\bibliographystyle{IEEEtran}
\bibliography{Ref}

\end{document}